\documentclass[]{amsart}
\usepackage{amssymb,amsmath}
\def\Fam{\mathcal{F}}
\def\ds{\displaystyle}
\def\ts{\textstyle}
\def\Wronskian#1#2#3{\left|\!\!\begin{array}{cc}
#1_1#2 &#1_2#2\\#1_1#3 & #1_2#3 \\
\end{array}\!\!\right|}

\begin{document}
\begin{center}{\Large\sc  Normal Families\\

\medskip

And Linear Differential Equations}\\

\bigskip

\sc by Norbert Steinmetz\\

\medskip

Technical University of Dortmund\end{center}

\bigskip

By Marty's Criterion (see Ahlfors \cite{A}), normality of any family
$\Fam$ of meromorphic functions on some domain $D$ is equivalent
to local boundedness of the corresponding family $\Fam^\#$ of {\it
spherical derivatives}
 $$f^\#=\frac{|f'|}{1+|f|^2}.$$
Recently, J.\ Grahl and S.\ Nevo(\footnote{I learned about this in
a talk given by J.\ Grahl at the second Bavarian-Qu\'ebec
Mathematical Meeting \& Tag der Funktionentheorie, November 22-27,
2010, University of W\"urzburg.}) proved a normality criterion
involving the spherical derivative by utilising the so-called
Zalcman Lemma (see L.\ Zalcman \cite{Z}). At first glance it looks
very surprising since it is based on a {\it lower} bound for the
spherical derivative.
\medskip

\noindent{\sc Theorem (Grahl \& Nevo \cite{GN}).} {\it Suppose all
functions of the family $\Fam$ satisfy $f^\#(z)\ge \epsilon$ for
some fixed $\epsilon>0.$ Then $\Fam$ is normal.}

\medskip The aim of this note is to give a
completely different proof, which has the advantage to yield
explicit upper bounds for $f^\#$. It is based on a property
equivalent to $f^\#(z)>0$, namely {\it local univalence} of the
function $f$, which again is equivalent to the fact that the
corresponding {\it Schwarzian derivative}
 $$S_f=\left(\frac{f''}{f'}\right)'-\frac 12\left(\frac{f''}{f'}\right)^2$$
is holomorphic on $D$.

\medskip
\noindent{\sc Theorem.} {\it Let $f$ be meromorphic on the unit
disc $\mathbb{D}$ satisfying $f^\#(z)\ge \epsilon>0.$ Then $f$ has
the form
 \begin{equation}\label{darstellung}f=\frac{w_1}{w_2},\end{equation}
where the functions $w_1$ and $w_2$ are holomorphic on
$\mathbb{D}$ and satisfy
\begin{equation}\label{beding}|w_1(z)|^2+|w_2(z)|^2 \le \frac1{\epsilon},
 \quad\Wronskian{w}{}{'}= 1,\quad{\rm and}\quad \Wronskian{w}{}{''}=0.\end{equation}
Moreover,
\begin{equation}\label{abschaetz}f^\#(z)
\le\frac{2/\epsilon}{(1-|z|)^2}\quad {\rm and}\quad|S_f(z)|
 \le\frac{4/\epsilon}{(1-|z|)^3} \end{equation}
hold on $\mathbb{D}$.}

\medskip\noindent{\sc Proof.} Since $f^\#$ is non-zero, $f$ is locally
univalent and its Schwarzian derivative is holomorphic on
$\mathbb{D}$. It is well known that this implies the
representation (\ref{darstellung}), where $w_1$ and $w_2$ form a
fundamental set of the linear differential equation
 \begin{equation}\label{diffeq}w''+{\ts\frac 12}S_f(z)w=0.\end{equation}
Then the third condition in (\ref{beding}) always holds
(reflecting the fact that the coefficient of $w'$ vanishes
identically), hence the Wronskian of any two solutions is
constant. To make some definite choice we normalise by the second
condition in (\ref{beding}), which makes the pair $(w_1,w_2)$
unique up to sign and from which
 $$f'=\frac{-1}{w_2^2} \quad{\rm and}\quad f^\#=\frac 1{|w_1|^2+|w_2|^2},$$
hence the first condition in (\ref{beding}) follows. To prove
(\ref{abschaetz}) we just remark that from $\Wronskian{w}{}{'}= 1$
and the Cauchy-Schwarz inequality follows
 $$f^\#=\frac1{|w_1|^2+|w_2|^2}\le |w'_1|^2+|w'_2|^2,$$
while
 ${\ts\frac 12}S_f=\Wronskian{w}{'}{''}$
yields
 $${\ts\frac 12}|S_f|\le |w'_1||w''_2|+|w_1''||w'_2|.$$
The standard Cauchy estimate
 $$|w|\le 1/\sqrt{\epsilon}\quad\Rightarrow\quad|w'|\le
\frac{1/\sqrt{\epsilon}}{1-|z|}\quad{\rm and}\quad |w''|\le
 \frac{1/\sqrt{\epsilon}}{(1-|z|)^2}$$
then gives the estimate in both cases of
(\ref{abschaetz}).\hfill$\square$

\bigskip\noindent{\sc Remarks and Questions.}\smallskip

$\bullet$ For $\epsilon>0$ fixed, the family $\Fam_\epsilon$ of
all functions $f$ satisfying $f^\#\ge\epsilon$, and also the
family $S_{\Fam_\epsilon}$ of corresponding Schwarzian derivatives
is compact, and
$$\Phi_\epsilon(r)=\sup\{f^\#(z): |z|\le r,~f\in\Fam_\epsilon\}
 \le 2\epsilon^{-1}(1-r)^{-2}\quad(0\le r<1)$$
holds. To obtain a lower bound for $\Phi_\epsilon$ we consider
$f(z)\ds=\Big(\frac{1+z}{1-z}\Big)^{i\lambda}$ (Hille's example
\cite{hille} showing that Nehari's univalence criterion
\cite{nehari} is sharp). It has spherical derivative
$f^\#(z)=\ds\frac{\lambda}{|1-z^2|}\frac2{|f(z)|+|f(z)|^{-1}}$ and
Schwarzian derivative $S_f(z)=2(1+\lambda^2)(1-|z|^2)^{-2}$, and
satisfies $f^\#(z)>f^\#(\pm i)=\lambda/\cosh \frac\pi2\lambda$ and
$f^\#(x)=\lambda(1-x^2)^{-1}$ $(-1<x<1),$ from which
$$\Phi_\epsilon(r)\ge \big(\log(1/\epsilon)+O(\log\log(1/\epsilon))\big)
 \;(1-r)^{-1}\quad(0<\epsilon<\epsilon_0,~0< r<1)$$
follows ($\epsilon_0\approx 0.42$ is the maximum of $\lambda/\cosh
\frac\pi2\lambda$ in $0<\lambda<\infty$). The true value of
$\Phi_\epsilon(r)$ has to remain open; is it
$C(\epsilon)(1-r)^{-1}$? The problem to determine
$\sup_{\Fam_\epsilon}|S_f(z)|$ also remains open.\smallskip

$\bullet$ Since $|w_1|^2+|w_2|^2$ is subharmonic, the spherical
derivative satisfies the mini\-mum principle: If $f$ is
meromorphic on some domain $D$, then $f^\#$ has no local minima
except at the critical points of $f$ (see also \cite{GN},
Prop.~4.) Actually, $-\log f^\#$ is subharmonic off the zeros of
$f^\#$.\smallskip

$\bullet$ By Thm.\ 3 of \cite{GN}, $\Fam_{\epsilon}=\emptyset$ if
$\epsilon>1/2$, while $\Fam_{1/2}$ consists of the rotations of
the Riemann sphere. Using $\Wronskian{w}{}{'}=1$ and
$1/f^\#=|w_1|^2+|w_2|^2\ge 2|w_1||w_2|$, this may be shown as
follows (similar to \cite{GN}): We first suppose $w_1(0)=0$ and
set $v(z)=w_2(z)w_1(z)/z$. Then $v(0)=w'_1(0)w_2(0)=-1$, hence
$\max\limits_{|z|=r}|v(z)|\ge 1$ holds by the maximum principle,
this implying $\min\limits_{|z|=r}|z|f^\#(z)\le 1/2$ and
$\inf\limits_\mathbb{D} f^\#(z)\le 1/2$. Also $\inf
\limits_\mathbb{D} f^\#(z)=1/2$ gives $|v(z)|\le 1=|v(0)|$, thus
$v(z)\equiv v(0)=-1$. This, however, is only possible if
$w_1(z)=cz$ and $w_2(z)=-1/c$, hence $f(z)=c^2z$, and from
$\inf\limits_\mathbb{D} f^\#(z)=|c|^2/(1+|c|^4)$ then follows
$|c|=1$. Without the normalisation $f(0)=0$,
$\inf\limits_\mathbb{D} f^\#=1/2$ implies that $f$ is a rigid
motion of the sphere.

\smallskip$\bullet$ The upper bound for $f^\#$
may be slightly improved. Given $z\in\mathbb{D}$ we may assume
$f(z)=0$ by rotating the Riemann sphere. We then have
$f^\#(z)=|w_2(z)|^{-2}$ and $w'_1(z)w_2(z)=-1$, hence
$f^\#(z)=|w'_1(z)|^{2}$. By the Schwarz-Pick lemma (thanks to J.\
Grahl for the keyword) applied to $\sqrt{\epsilon}w_1$ we thus
obtain
 $$f^\#(z)=|w'_1(z)|^{2}\le\frac{1/\epsilon}{(1-|z|^2)^2}.$$

\smallskip$\bullet$ The representation (\ref{darstellung}) together
with the first condition in (\ref{beding}) implies that $f$ has
bounded {\it Nevanlinna characteristic} ($f$ is called {\it of
bounded type}), so that by the {\it Ahlfors-Shimizu formula}
$$\lim_{r\to 1}T(r,f)=\frac1\pi\int_0^1(1-\rho)\int_0^{2\pi}f^\#(\rho
 e^{i\theta})^2\,d\theta\,d\rho$$
is finite (see Nevanlinna \cite{N}). We note, however, that there
are functions of bounded type having spherical derivative growing
arbitrarily fast, see \cite{aul}. Thus, although normal functions
[satisfying  $f^\#(z)=O((1-|z|)^{-1})$ as $|z|\to 1$] have
Nevanlinna characteristic $T(r,f)=O(-\log (1-r))$, there are
functions of bounded type that are not normal.

\noindent{\small Adress: Institut f\"ur Mathematik, D-44221 Dortmund, Vogelpothsweg 87, Germany\\
E-mail: {stein@math.tu-dortmund.de\\
Homepage: www.mathematik.tu-dortmund.de/steinmetz/}}

\begin{thebibliography}{99}
\bibitem{A} L.\ V.\ Ahlfors, {\it Complex Analysis},
McGraw-Hill 1979.
\bibitem{aul} R.\ Aulaskari and J. R\"atty\"a, {Nevanlinna class contains functions whose
sperical derivatives grow arbitrarily fast}, {\it Ann.\ Acad.\
Sci.\ Fenn.\ Math.} {\bf 34}, 387 - 390 (2009).
\bibitem{GN} J.\ Grahl and S.\ Nevo, A note on spherical derivatives and normal families,
{\it arXiv: 1010.4654v1} [math.CV], 12 p.\ (2010).
\bibitem{hille} E.\ Hille, {Remarks on a paper by Zeev Nehari,}
{\it Bull.\ Amer.\ Math.\ Soc.} {\bf 55}, 552 - 553 (1949).
\bibitem{nehari} Z.\ Nehari, {The Schwarzian derivative and schlicht functions,}
{\it Bull.\ Amer.\ Math.\ Soc.} {\bf 55}, 544-551 (1949).
\bibitem{N} R.\ Nevanlinna, {\it Eindeutige analytische
Funktionen,} Springer 1936.
\bibitem{Z} L.\ Zalcman, {A heuristic principle in
function
 theory}, {\it Amer. Math. Monthly} {\bf 82}, 813-817 (1975).
\end{thebibliography}
\end{document}